\documentclass[preprint]{imsart}

\usepackage{color}
\usepackage{amsthm,amsmath, amssymb}
\usepackage[authoryear]{natbib}
\RequirePackage{graphicx}

\startlocaldefs

\newcommand{\mathbold}[1]{\mbox{\boldmath $#1$}}

\endlocaldefs

\begin{document}

\begin{frontmatter}

\title{Quantifying the cost of simultaneous non-parametric
  approximation of several samples}
\runtitle{Quantifying the cost of simultaneous approximation}


\author{\fnms{P.~L.} \snm{Davies}\ead[label=e1]{laurie.davies@uni-due.de}\thanksref{t1}},
\address{University of Duisburg-Essen\\Technical University Eindhoven\\\printead{e1}} 
\author{\fnms{A.} \snm{Kovac}\ead[label=e2]{A.Kovac@bristol.ac.uk}\thanksref{t1}},
\address{University of Bristol\\\printead{e2}}
\thankstext{t1}{Research supported in part by Sonderforschungsbereich 475, Technical University of Dortmund}
\affiliation{University Duisburg-Essen, Technical University of
  Eindhoven and University of Bristol}

\runauthor{Davies and Kovac}

\begin{abstract}
We consider the standard non-parametric regression model with Gaussian
errors but where the data consist of different samples. The
question to be answered is whether the samples can be adequately
represented by the same regression function. To do this we define for
each sample a universal, honest and non-asymptotic confidence region
for the regression function. Any subset of the samples can be
represented by the same function if and only if the intersection of
the corresponding confidence regions is non-empty. If the
empirical supports of the samples are disjoint then the
intersection of the confidence regions is always non--empty and a
negative answer can only be obtained by placing shape or quantitative
smoothness conditions on the joint approximation. Alternatively a
simplest joint approximation function can be calculated which gives a
measure of the cost of the joint approximation, for example, the
number of extra peaks required.
\end{abstract}



\end{frontmatter}

\newcounter{fig1}
\newcounter{fig2}
\section{Introduction} \label{introduction}
We consider the following problem in non-parametric regression. Given
$k$ samples     
\begin{equation} \label{samples}
{\mathbold y_{in_i}}=\{(t_{ij},y_{ij}):j=1,\ldots,n_i\},\,
i=1,\ldots,k,
\end{equation}
with supports 
\begin{equation}\label{support}
S_{in_i}=\{t_{i1}<t_{i2}<\ldots <t_{in_i}\},\,i=1,\ldots,k,
\end{equation}
the question to be answered is whether they can be simultaneously
represented by a common function $f.$ The standard approach is to 
assume that the data sets were generated according to the model
\begin{equation} \label{basicmod}
Y_{in_i}(t)=f_i(t)+\sigma_i Z_i(t),\,\, i=1,\ldots,k,\quad t \,\in\,
[0,\,1], 
\end{equation}
and then to consider the null and alternative hypotheses   
\begin{equation} \label{H0}
H_0: f_1=\ldots =f_k \quad H_1: f_i\ne f_j\quad \text{for some } i,j.
\end{equation}
We assume that the noise processes $Z_i(t), i=1,\ldots,k$
are independent and  standard Gaussian white noise. Individual samples
generated under (\ref{basicmod}) will be denoted by   
\[{\mathbold Y_{in_i}}=\{(t_{ij},Y_{ij}):j=1,\ldots,n_i\},\,
i=1,\ldots,k.\]
Here and in the following we use minuscule letters to denote
general data sets and majuscule letters for data generated under
(\ref{basicmod}).  We shall mostly restrict attention to the case 
$k=2;$ the extension to more samples poses no problems. 

Within this setup it is possible to construct tests which are
asymptotically consistent if $\lim n_i=\infty,\, i=1,2$, and which can
detect alternatives  converging to the null hypothesis at certain
rates. This may be formalized by putting
\begin{equation}\label{diffn}
f_1(t)-f_2(t)=f_{1,n_1}(t)-f_{2,n_2}(t)=\Delta_n(t), \quad n =\min
(n_1,n_2)
\end{equation}
where $\Delta_n$ is a difference function and measures the rate of
convergence to the null hypothesis. The best result seems to be that of
\cite{NEUDET03} who  construct a test which can detect alternatives
which converge to the null hypothesis at the optimal rate
$\Delta_n=\text{O}(n^{-1/2}).$  If the supports are equal,
$S_{in_i}=\{t_1,\ldots,t_{n_i}\},\, i=1,2,$ then 
it is not difficult to construct such a test as the differences 
$Y_{1n_1}(t_j)-Y_{2n_2}(t_j)$ do not depend on $f$ (see for example
\cite{DEL92} and \cite{FANLIN98}). The result of \cite{NEUDET03}
continues to hold even if the supports are disjoint, $S_{1n_1}\cap
S_{2n_2}=\emptyset.$  In this case, however, there are difficulties
which can be most clearly seen in the case of exact data    
\[y_{ij}=f_i(t_{ij}), \,\, t_{ij} \in S_{in_i},\quad i=1,2.\] 
If we denote the supremum norm on $[0,\,1]$ by $\Vert\cdot
\Vert_{\infty}$ then the null and alternative hypotheses of (\ref{H0})
may be rewritten as    
\begin{equation} \label{H0ref}
H_0: \Vert f_1-f_2\Vert_{\infty}=0, \quad H_1:\Vert
f_1-f_2\Vert_{\infty} >0.  
\end{equation}
If the values of $f_1$ and $f_2$ are known only on disjoint sets
$S_{1n_1}$ and $S_{2n_2}$ respectively,  then it is not possible to
decide between $H_0$ and $H_1$. This continues to hold even if $f_1$ 
and $f_2$ are subject to qualitative smoothness conditions such as
infinite differentiability: all one does is to interpolate the data
points using such a function. The addition of noise and the use of
asymptotics does not solve the problem as indicated by Figure
\ref{fig1}. The top panel shows two data sets of sizes $n_1=n_2=500$
generated according to $Y_1(t)=\exp(1.5t)+0.25Z_1(t)$ and 
$Y_2(t)=\exp(1.5t)+3+0.25Z_2(t)$ with disjoint supports taken to be
i.i.d. uniform random variables on $[0,1].$  The centre panel shows a
joint piecewise constant approximating function with 514 local extreme
values. It can be made infinitely differentiable by convolving it with
a Gaussian kernel with a small bandwidth. The bottom panel shows a
sample of size $n=1000$ generated using the function of the centre
panel. It looks very much like the two original data sets.

In order to distinguish between $H_0$ and $H_1$ it is necessary to
place either quantitative conditions on $f_1$ and  $f_2$ such as $\Vert
f_1^{(1)}\Vert_{\infty} \le 1,\,\Vert f_2^{(1)}\Vert_{\infty} \le 1,$
or shape restrictions such as $f_1$ and $f_2$ being monotone.
In spite of this all conditions imposed in the
literature are of a qualitative form: \cite{HALLHART90}, a bounded
first derivative; \cite{HAEMAR90}, H\"older continuity;
\cite{KINHARWEH90}, `at least uniform continuity´; \cite{KUL95},
\cite{KUSWANG97}, a continuous second derivative; \cite{MUNDETT98},
H\"older continuity of order $\beta > 1/2$; \cite{DETTNEUM01}, a
continuous $r$th derivative: \cite{LAV01}, a second derivative which
is uniformly Lipschitz of order $\beta,\, 0 \le \beta \, <1$;
\cite{NEUDET03}, continuous derivatives of order $d \ge 2$. The
problem is one of uniform convergence which is required to make the
results applicable for finite $n$ and which does not follow from
qualitative conditions alone. What can be said is that if the functions
differ, then any joint approximation will become more complicated as
the sample sizes increase. It is this increase in complexity which we
call the cost of the simultaneous approximation. This is shown in
Figure 1 where the individual approximations are monotone (top panel)
but the simplest joint approximation has 514 local extreme values
(centre panel). In the remainder of the paper we show how the
quantification can be carried out. Our approach can be split into two
parts:

\begin{itemize}
\item[(1)]Firstly, for each sample ${\mathbold
  y_{in_i}}$ we define a so called approximation region ${\mathcal
  A}_{in_i}$ which specifies those functions $f_i$ for which the model
(\ref{basicmod}) is an adequate approximation for the sample. The
intersection of the approximation regions ${\mathcal
  A}_{1,n_1}\cap{\mathcal A}_{2,n_2} $ contains all those functions
which simultaneously approximate both samples. It is also the
approximation region for the simultaneous approximation. A similar
idea in the context of the one-way table in the analysis of variance
is expounded in \cite{DAV04}.  
\item[(2)] Secondly, using some measure of complexity we regularize
  within each approximation region by choosing the simplest function
  which is consistent with the data. This is in the spirit of
  \cite{DON88} who pointed out that in non-parametric regression and
  density problems it is possible only to give lower bounds on certain
  quantities of interest such the number of modal values. 
\end{itemize}

In Section \ref{secapproxreg} we define the approximation or
confidence regions ${\mathcal A}_{in_i}$ and in Section \ref{secreg}
we apply the ideas and concepts to the problem of comparing regression
functions.

\section{Approximation regions} \label{secapproxreg}
\subsection{Single samples}
The following is based on \cite{DAVKOVMEI08}. We consider a single
sample of data ${\mathbold Y}_n=(t_i,Y(t_i))_1^n$ generated under the model 
\begin{equation} \label{basicmod1}
Y(t)=f(t)+\sigma Z(t)
\end{equation}
where we take the $t_i$ to be ordered. Based on this model we consider
two different approximation or confidence regions ${\mathcal A}_n$ and
${\mathcal A}_n^*$ defined as follows. For any function $g$ and 
any interval $I \subset [0,\,1]$ we put 
\begin{equation} \label{defw}
w(g,{\mathbold Y}_n,I)= \frac{1}{\sqrt{\vert I\vert}}\,\sum_{t_i \in
  I}(Y(t_i)-g(t_i)) 
\end{equation}
where $\vert I\vert$ denotes the number of points $t_i \in I.$ The
confidence region  ${\mathcal A}_n$ is defined by
\begin{equation} \label{confreg}
{\mathcal A}_n({\mathbold Y}_n, {\mathcal I}_n ,\sigma, \tau_n)= \{g:\,
\max_{ I \in {\mathcal I}_n} \, \vert w({g, \mathbold 
  Y}_n,I)\vert\ \le \sigma \sqrt{\tau_n \log(n)\,}\,\}.
\end{equation} 
where ${\mathcal I}_n$ is a collection of intervals of $[0,\,1]$.  We
restrict attention to the cases where ${\mathcal I}_n$ is either the
set of all intervals or a set of intervals of the form
\begin{eqnarray} 
\lefteqn{\hspace*{-1.5cm}{\mathcal
    I}_n(\lambda)=\left\{[t_{l(j,k)},\,t_{u(j,k)}]\,: \,l(j,k)=\lfloor 
(j-1)\lambda^k+1\rfloor,\right.} \nonumber\\
&&\left.u(j,k)=\min\{\lfloor j\lambda^k\rfloor,n\},j=1,\ldots,\lceil
n\lambda^{-k}\rceil,k=1,\ldots,\lceil \log n/\log
\lambda\rceil\right\}.\label{multischeme} 
\end{eqnarray}
for some $\lambda >1.$ Our default choice is the (wavelet) dyadic
scheme ${\mathcal I}_n(2)$. For any given $\alpha$ and collection of
intervals ${\mathcal I}_n$ 
we define $\tau_n(\alpha)$ by
\begin{equation} \label{taualpha}
{\mathbold P}\left( \max_{ I \in {\mathcal I}_n} \frac{1}{\sqrt{\vert
      I\vert}} \Big \vert \sum_{i \in I}Z(t_i)\Big\vert \le
  \sqrt{\tau_n(\alpha)\log n\,} \right) =\alpha.
\end{equation}
The value of $\tau_n(\alpha)$ may be determined by simulations.  These
show that for ${\mathcal I}_n={\mathcal I}_n(2)$ we have $\tau_n(0.95)
\le 3$ for all $n\ge 500.$ If ${\mathcal I}_n$ contains all singletons
$\{t_i\}$, as will always be the case, it follows from \cite{DUEMSPO01}
and \cite{KAB07} that $\lim_{n   \rightarrow \infty}\tau_n(\alpha)=2$
for any $\alpha.$   One immediate consequence of (\ref{taualpha}) is
\begin{equation} \label{confreg0}
{\mathbold P}( f \in {\mathcal A}_n({\mathbold Y}_n, {\mathcal I}_n,
\sigma, \tau_n(\alpha)))= \alpha
\end{equation}
so that ${\mathcal A}_n$ is a universal, exact and non-asymptotic
confidence region  for $f$ of size $\alpha.$

The confidence region  (\ref{confreg}) treats all intervals
equally. The second confidence region ${\mathcal A}_n^*$ downweights
the importance of small intervals and is defined as
follows. \cite{DUEMSPO01} extended L\`{e}vy's uniform modulus of
continuity of the Brownian motion and showed that
\begin{equation}\label{dumspok}
\sup_{0 < s < t <1} \frac{\frac{(B(t)-B(s))^2}{t-s}
  -2\log(1/(t-s))}{\log(\log(e^e/(t-s)))} < \infty \quad \text{a.s.}
\end{equation}
If we embed the partial sums $\sum_{i\in I}^j Z(t_i)/\sqrt{\vert I\vert},
\, I \in {\mathcal I}_n,$ in a standard Brownian motion it follows that
\begin{equation}\label{dumspokint}
\sup_{I \in {\mathcal I}_n} \frac{(\sum_{t_j\in I}Z(t_j))^2/\vert I\vert
-2\log(n/\vert I\vert)}{\log(\log(e^en/\vert I\vert)))}=\Gamma <
\infty \quad \text{a.s.}. 
\end{equation}
This implies that for any $\alpha$ we can find a
$\gamma_n=\gamma_n(\alpha)$ such that
\begin{eqnarray} 
\lefteqn{{\mathcal A}^*_n({\mathbold Y}_n, {\mathcal I}_n ,\sigma,
  \gamma_n(\alpha)) = \{g:\, \vert w({g, \mathbold
    Y}_n,I)\vert}\hspace{1cm} \label{confregbm} \\ 
&\le& \sigma \sqrt{2\log(n/\vert
   I\vert)+\gamma_n(\alpha)\log(\log(e^en/\vert I\vert))} \,\,
 \text{for all}\,\,I \in {\mathcal I}_n )\,\}.\nonumber
\end{eqnarray} 
is a universal, exact and non-asymptotic $\alpha$-confidence region for $f.$
The values of $\gamma_n$ may be determined by simulation. For
$\alpha=0.95$ and with ${\mathcal I}_n={\mathcal I}_n(2)$ a good
approximation for $\gamma_n(\alpha)$ for $n \ge 100$ is given by
\begin{equation} \label{gammaapprox}
\gamma_n(0.95)\approx 5.77-\exp(2.89-0.6\log(n)).
\end{equation}

The confidence regions ${\mathcal A}_n({\mathbold Y}_n, {\mathcal
  I}_n, \sigma, \tau_n)$ and ${\mathcal A}^*_n({\mathbold Y}_n,
{\mathcal I}_n ,\sigma, \gamma_n)$ both require the true value of
$\sigma.$ We indicate how this may be obtained from the data in such a
manner that the confidence region now becomes honest (\cite{LI89})
rather than exact. The following argument makes the somewhat casual
remarks on the problem made in \cite{DAVKOVMEI08} more precise. Using
the normal approximation for the binomial $(n,1/2)$ distribution it
follows that for an i.i.d. sample $(W_1,\ldots,W_n)$ with common
continuous distribution $P$ with median  $\text{med}(P)$
\[ {\mathbold P}\big( W_{(\lceil n/2+z_{\beta}\sqrt{n}/2\rceil)}\ge
  \text{med}(P)\big)=\beta\]
where $W_{(i)}$ denotes the $i$th order statistic of the sample and
$z_{\beta}$ the $\beta$-quantile of the standard normal
distribution. On putting $\beta=0.995$ we obtain
\[ {\mathbold P}\big( W_{(\lceil n/2+1.288\sqrt{n}\,\rceil)}\ge
  \text{med}(P)\big)=0.995.\]
We now apply this to the $\lfloor n/2\rfloor$ random variables 
\[V_i=\vert Y_{2i-1}-Y_{2i}\vert,\quad i=1,\ldots, \lfloor n/2\rfloor.\]
It follows from  \cite{ANDER55} that whatever the value of the
function $f$ 
\[{\mathbold P}\big( V_i\ge 1.4826\sigma/\sqrt{2}\Big) \ge 1/2\]  
and consequently
\[ {\mathbold P}\big( V_{(\lceil n/4+0.9108\sqrt{n}\,\rceil)}\ge
 1.4826\sigma/\sqrt{2}\big)=0.995.\]
On using the corresponding result for the  $\lfloor (n-1)/2\rfloor$
random variables 
\[\vert Y_{2i}-Y_{2i+1}\vert,\quad i=1,\ldots, \lfloor (n-1)/2\rfloor
\]
it follows that if we define ${\hat \sigma_n}$ to be the $\lceil
n/2+1.814\sqrt{n}\,\rceil$ order statistic of the random variables
\[ \frac{1.4826}{\sqrt{2}}\vert
Y(t_2)-Y(t_1)\vert,\ldots,\frac{1.4826}{\sqrt{2}}\vert
Y(t_n)-Y(t_{n-1})\vert \]  
then 
\[ {\mathbold P}\big({\hat \sigma_n}\ge \sigma\big) \ge 0.99\]
for all $n \ge 100$ say whatever the value of $f$. It follows that 
 ${\mathcal A}_n({\mathbold Y}_n, {\mathcal
  I}_n, {\hat \sigma_n}, \tau_n)$ and ${\mathcal A}^*_n({\mathbold Y}_n,
{\mathcal I}_n ,{\hat \sigma_n}, \gamma_n)$ are now universal and
non-asymptotic honest confidence regions whatever the value of $f$ but
with $\alpha$ replaced by $\alpha -0.01,$
\begin{equation} \label{confreg0hon}
{\mathbold P}\big( f \in {\mathcal A}_n({\mathbold Y}_n, {\mathcal I}_n,
{\hat \sigma}_n, \tau_n(\alpha))\big)\ge \alpha-0.01.
\end{equation}
with the corresponding inequality for ${\mathcal A}_n^*.$ In spite of
this the default value for ${\hat \sigma}_n$ we shall use in this
paper is
\begin{equation} \label{sigma}
{\hat \sigma}_n=  \frac{1.4826}{\sqrt{2}}\text{median}\,(\vert
Y(t_2)-Y(t_1)\vert,\ldots,\vert Y(t_n)-Y(t_{n-1})\vert).
\end{equation}
It is simpler, the difference is in general small, it was used in
\cite{DAVKOV01}, \cite{DAVGATWEI06}, \cite{DAGAMEMEMI08}
and it also corresponds to using the first order Haar wavelets  
to estimate $\sigma$ (\cite{DONJOHKERPIC95}).  

In \cite{DAV95} implicit use is made of an confidence
region based on the lengths of runs of the signs of the
residuals. Explicit universal, honest and non-asymptotic confidence
regions based on the signs of the residuals are to be found in
\cite{DUEM98A,DUEM03, DUEM07} and \cite{DUEMJO04}.

\subsection{A one-way table for regression functions}
This section extends the approach given in \cite{DAV04} for the
one-way table to the case of regression functions. We consider $k$
samples ${\mathbold Y}_{in_i}=(t_{ij},Y_i(t_{ij}))_{j=1}^{n_i}$
generated under (\ref{basicmod}). As a first step we replace the
$\alpha$ in (\ref{taualpha}) and (\ref{confregbm}) by
$\alpha_k=\alpha^{1/k}$ where $k$ is the number of samples. This
adjusts the size of each confidence region to take into account the
number of samples. The confidence region for the $i$th sample is given
by \begin{eqnarray} 
\lefteqn{{\mathcal A}_{in_i}={\mathcal A}_{in_i}({\mathbold Y}_{in_i},
  {\mathcal I}_{in_i} ,{\hat \sigma}_{in_i},\tau_{in_i}(\alpha_k))=}
\label{confregi}\hspace{2cm}\\ 
&&\big\{g:\, \max_{ I \in {\mathcal I}_{in_i}}
  \, \vert w({g, \mathbold Y}_{in_i},I)\vert\ \le {\hat \sigma}_{in_i}
  \sqrt{\tau_{in_i}(\alpha_k) \log(n_i)\,}\,\,\big\}. \nonumber
\end{eqnarray} 
We denote by ${\mathbold P}_{\mathbold f}$ with ${\mathbold
f}=(f_1,\ldots,f_k)$ the probability model where all the samples
${\mathbold Y}_{in_i},\, i=1,\ldots,k,$ are independently distributed
and ${\mathbold Y}_{in_i}$ was generated under (\ref{basicmod}) with
$f=f_i,\, i=1,\ldots,k$. It follows from the choice
$\alpha_k=\alpha^{1/k}$ that  
\begin{equation} \label{confreg3}
{\mathbold P}_{\mathbold f}\,( f_i \in {\mathcal A}_{in_i}({\mathbold
  Y}_{in_i}, {\mathcal I}_{in_i}
,{\hat \sigma}_{in_i},\tau_{in_i}(\alpha_k)),\, i=1,\ldots k)\ge \alpha \quad
\text{for all}\quad {\mathbold f}.
\end{equation}
All questions concerning the relationships between the functions $f_i$
can now be answered by using the confidence regions ${\mathcal
  A}_{in_i}$. For example, the question as to whether the $f_i$ are
all equal translates into the question as to whether 
\begin{equation} \label{jointconf}
{\mathcal A}_{\mathbold n_k}=\cap_{i=1}^k {\mathcal A}_{in_i}=
\cap_{i=1}^k {\mathcal A}_{in_i}({\mathbold Y}_{in_i}, {\mathcal
  I}_{in_i} ,{\hat \sigma}_{in_i}, \tau_{in_i}(\alpha_k)), \quad {\mathbold
  n_k}=(n_1,\ldots,n_k) 
\end{equation}
is empty or not. If the supports $S_{in_i}$ of the samples are not
disjoint then it is possible that the linear inequalities which define
the confidence regions are inconsistent. In this case ${\mathcal
  A}_{\mathbold n_k}=\emptyset$ and there is no joint approximating
function. If the supports $S_{in_i}$ of the samples are pairwise
disjoint then ${\mathcal A}_{\mathbold n_k}$ is non--empty and so there
always is a joint approximation function. Without further restrictions
on the joint approximating function nothing more can be
said. If however the joint approximating function is required to
satisfy, for example, a shape constraint such as monotonicity, then it
may be the case that there is no joint approximating function. Figure
\ref{fig1} shows just such a case where there are monotonic
approximations for each sample individually but no monotonic joint
approximation. To answer questions of this nature we must regularize
within ${\mathcal A}_{\mathbold n_k}$ and this is the topic of the
next section.

\section{Regularization}\label{secreg}
\subsection{Disjoint supports}
We consider firstly the case when the supports $S_{in_i},
i=1,\ldots, k$, are pairwise disjoint. In this case the joint
approximation region ${\mathcal A}_{\mathbold n_k}$ is non-empty and
will in general include many functions which would not be regarded as
being acceptable. Indeed, it may be that ${\mathcal A}_{\mathbold
  n_k}$ does not contain any acceptable function. The definition of
`acceptable' will usually be formulated in terms of shape or
quantitative smoothness constraints. 

\begin{figure}[h]
\begin{center}
\includegraphics[height=5cm,width=15cm,angle=0]{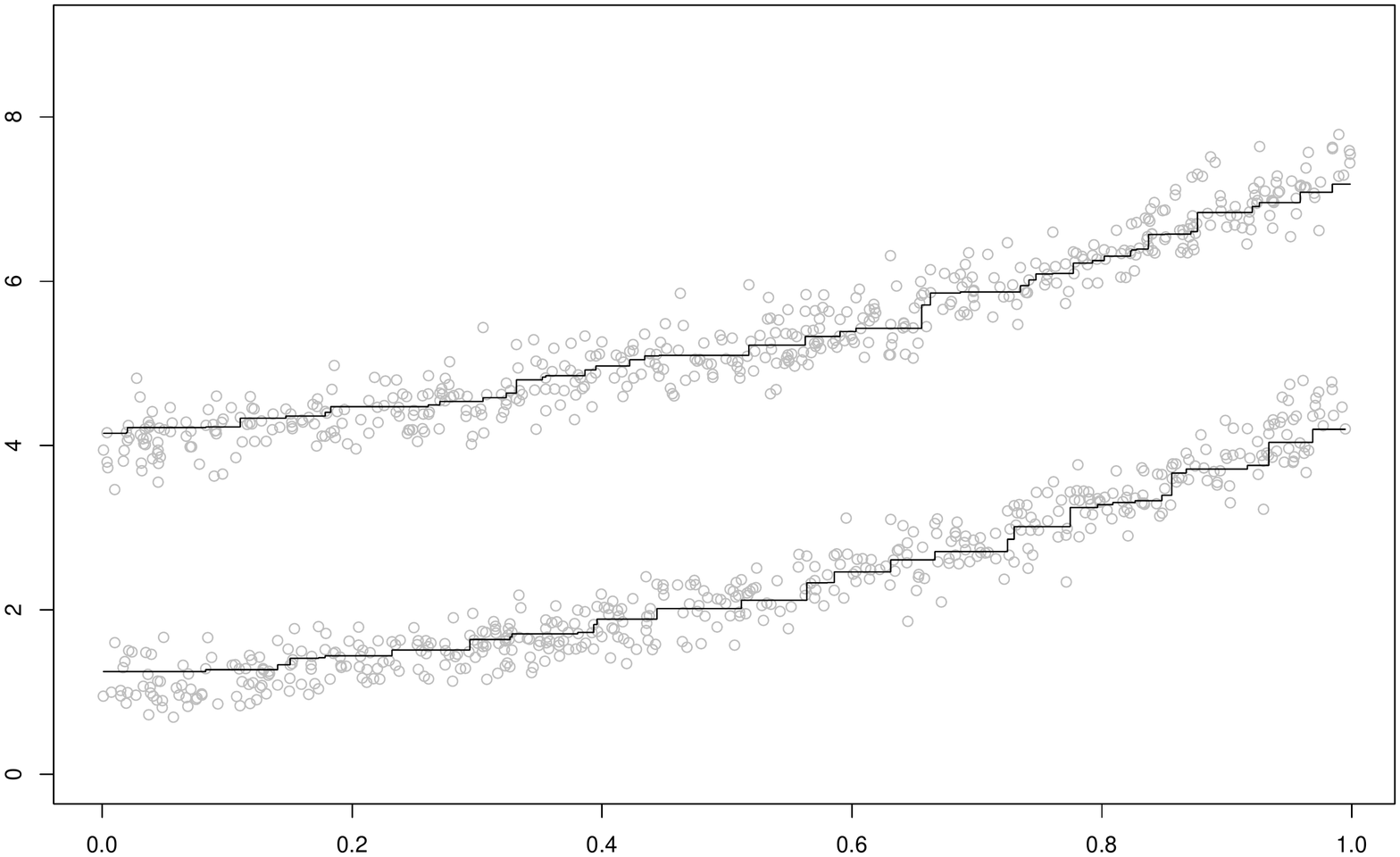}
\includegraphics[height=5cm,width=15cm,angle=0]{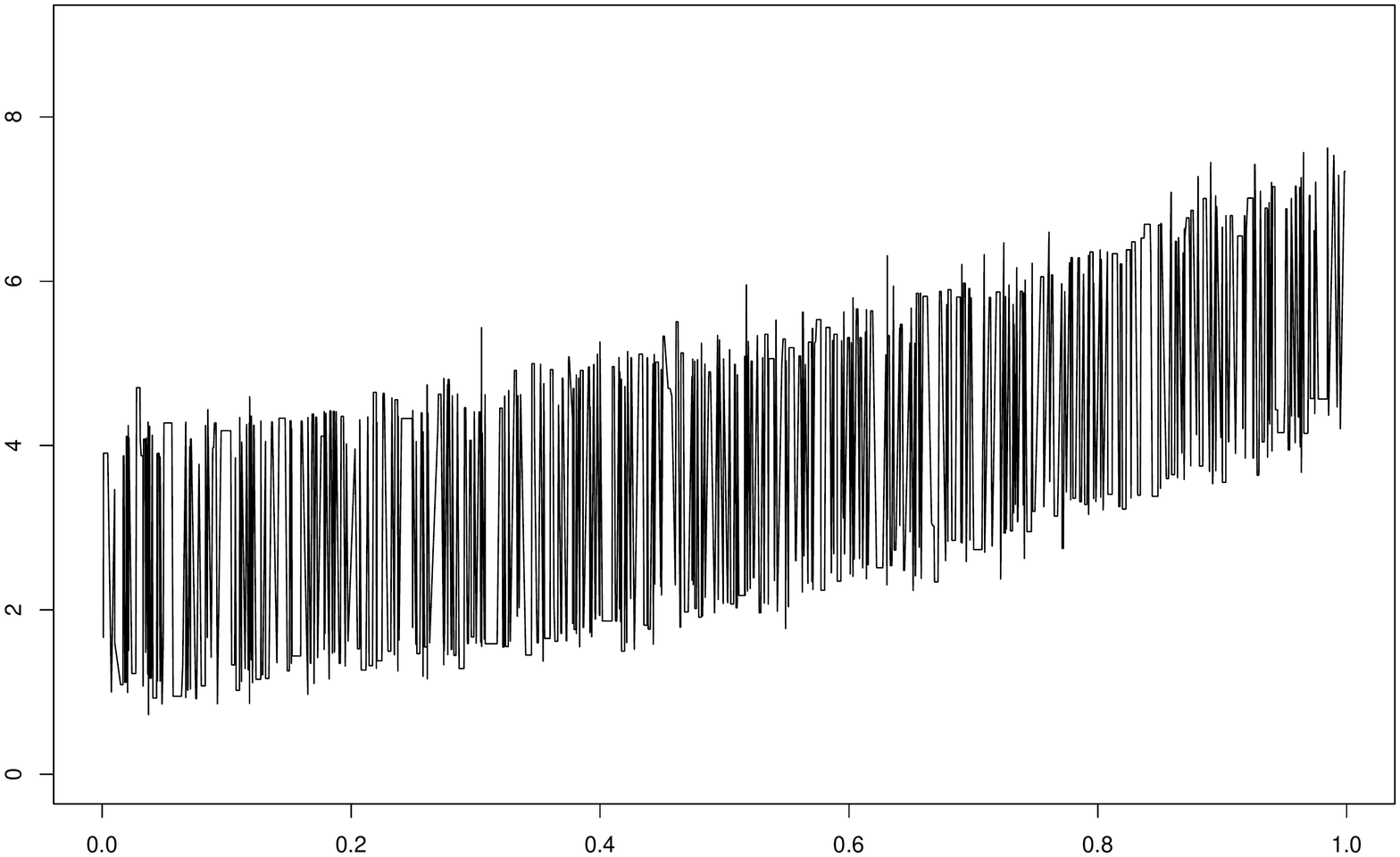}
\includegraphics[height=5cm,width=15cm,angle=0]{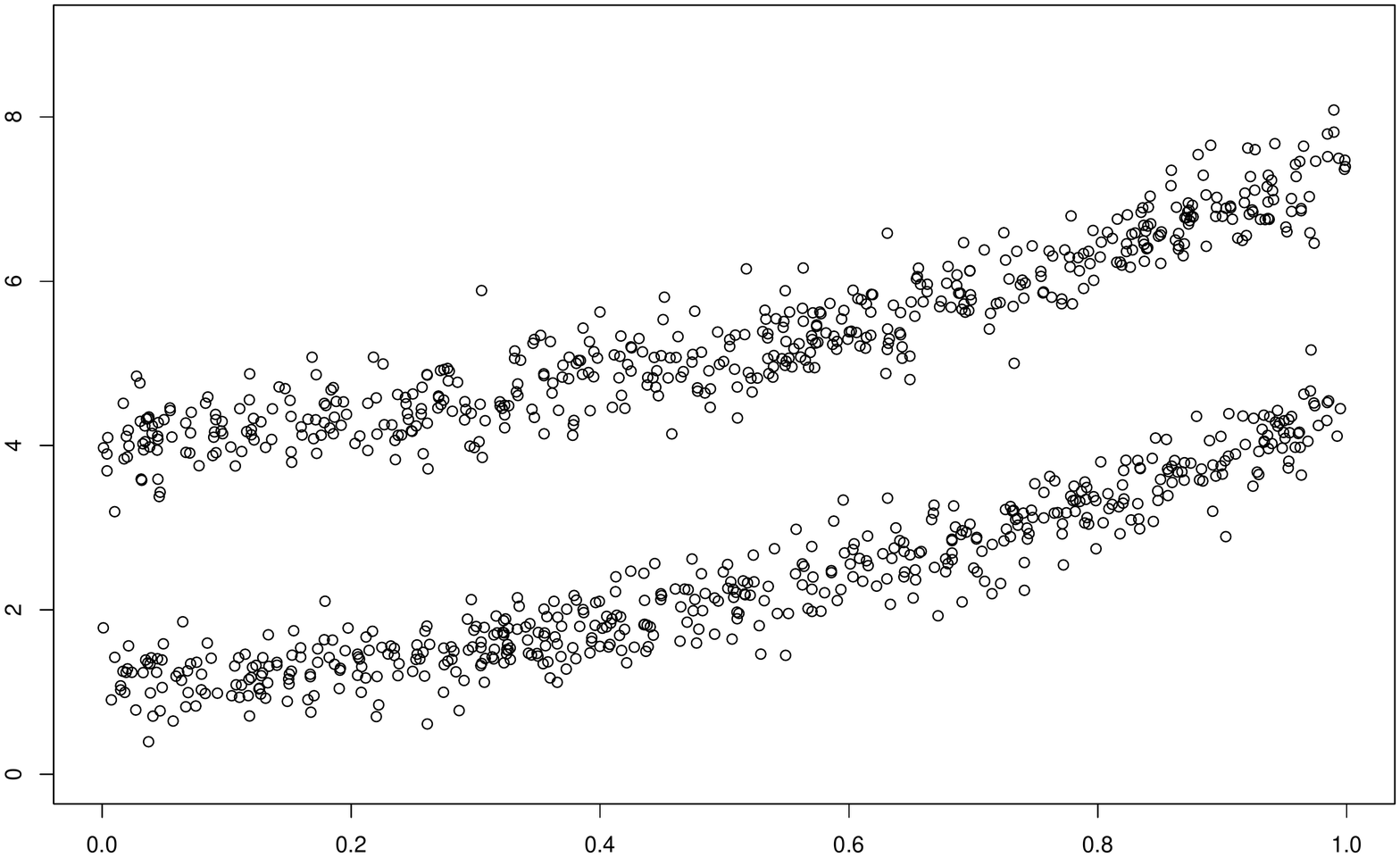}
\end{center}
\caption{The top panel shows two samples each of size 500 generated by
  $Y_1(t)=\exp(1.5t)+0.25Z(t)$ and $Y_2(t)=\exp(1.5t)+3+0.25Z(t)$
  together with the approximating monotonic curves. The design points
  were taken to be i.i.d. random variables uniformly distributed on
  $[0,\,1].$ The centre panel shows a joint approximating function with
  514 local extreme values. The bottom panel shows a sample of size
  $n=1000$ generated using the function of the centre panel. \label{fig1}}
\end{figure}

\begin{figure}[hb]
\begin{center}
\includegraphics[height=8cm,width=15cm,angle=0]{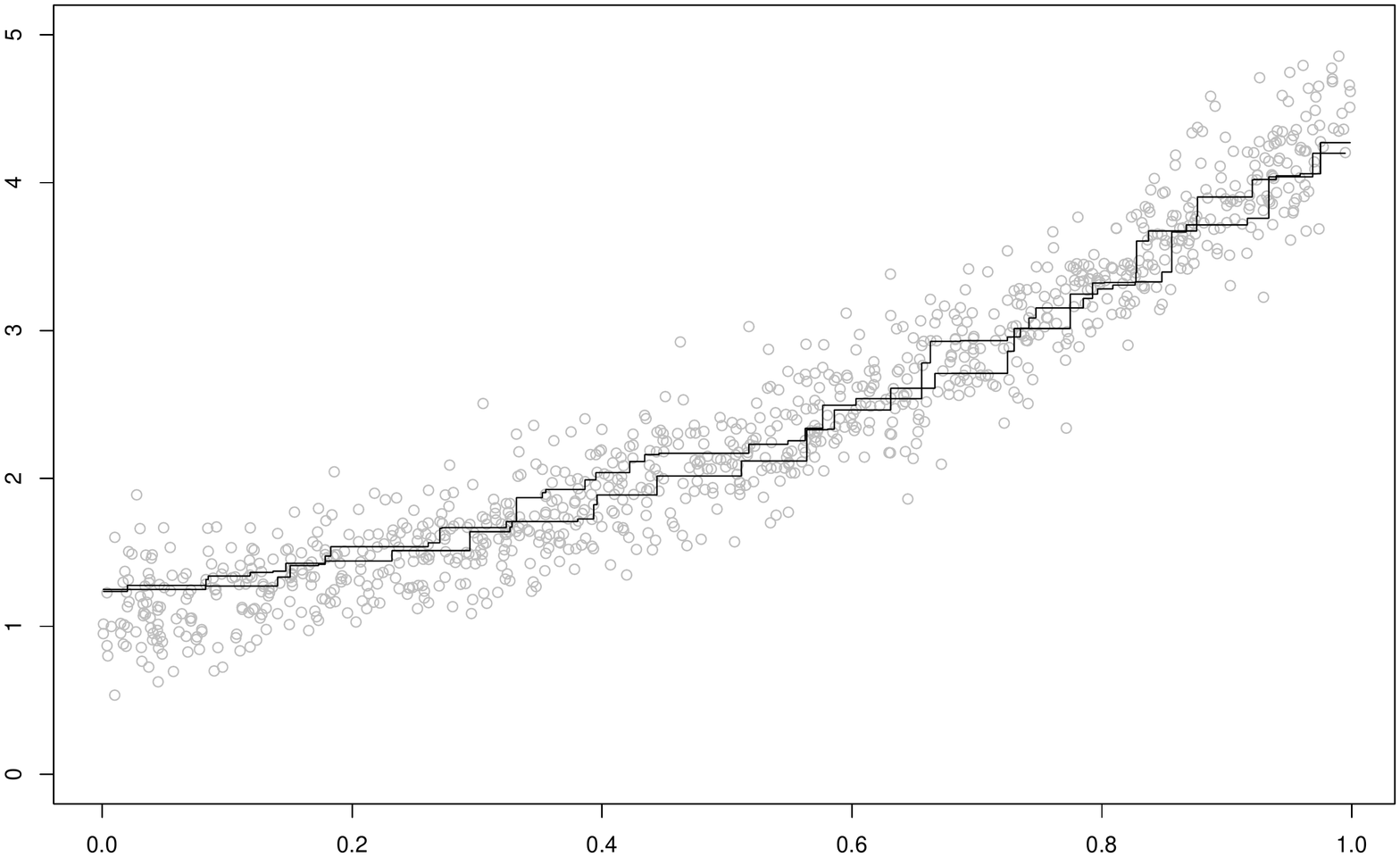}
\includegraphics[height=8cm,width=15cm,angle=0]{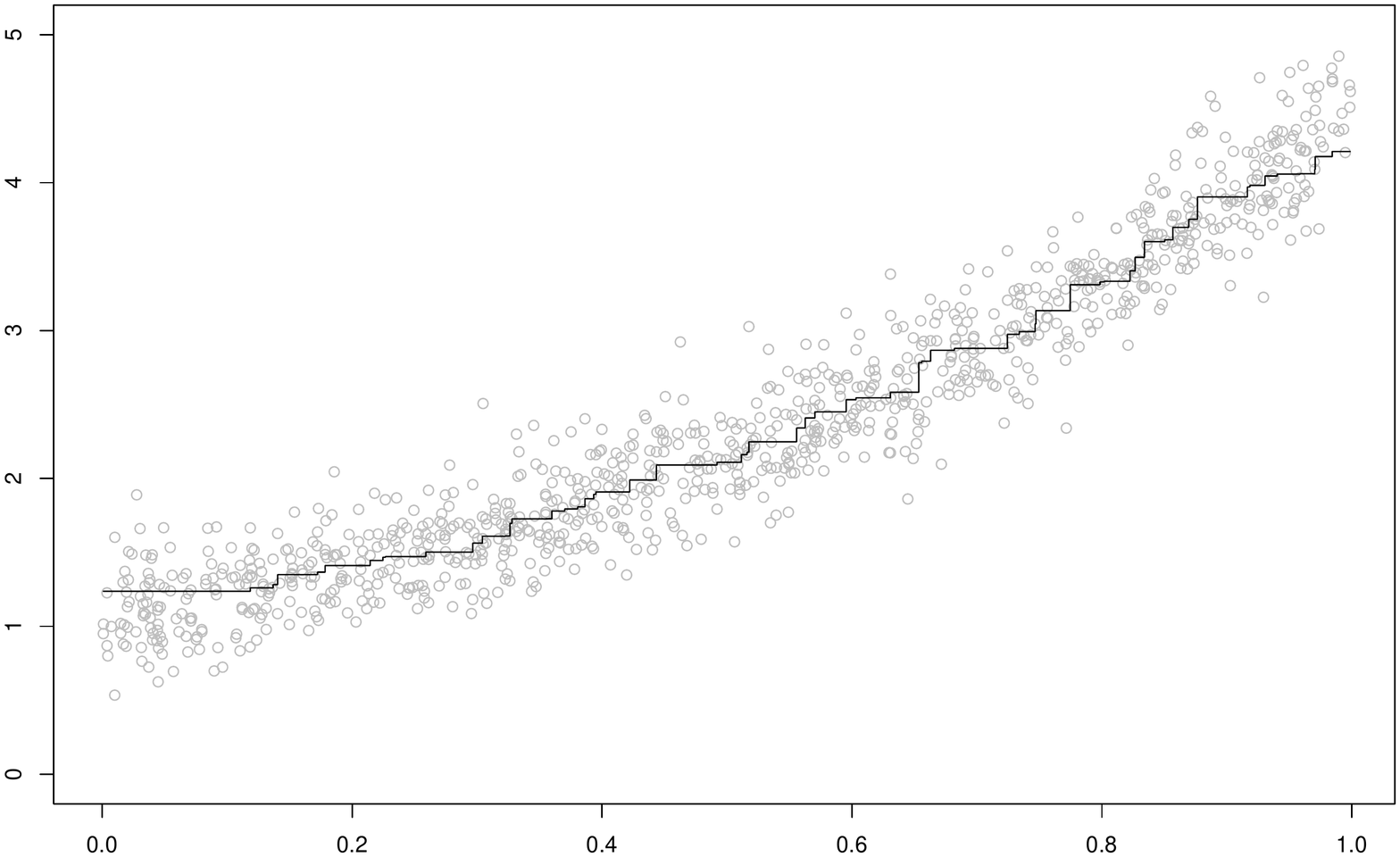}
\end{center}
\caption{The upper panel shows the data of Figure \ref{fig1} but with the
  values of $Y_2$ now given by $Y_2(t)=\exp(1.5t)+0.07+0.25Z(t).$
  There is now a joint monotonic approximating function which is shown
  in the lower panel. \label{fig2}} 
\end{figure}

\begin{figure}[hb]
\begin{center}
\includegraphics[height=8cm,width=15cm,angle=0]{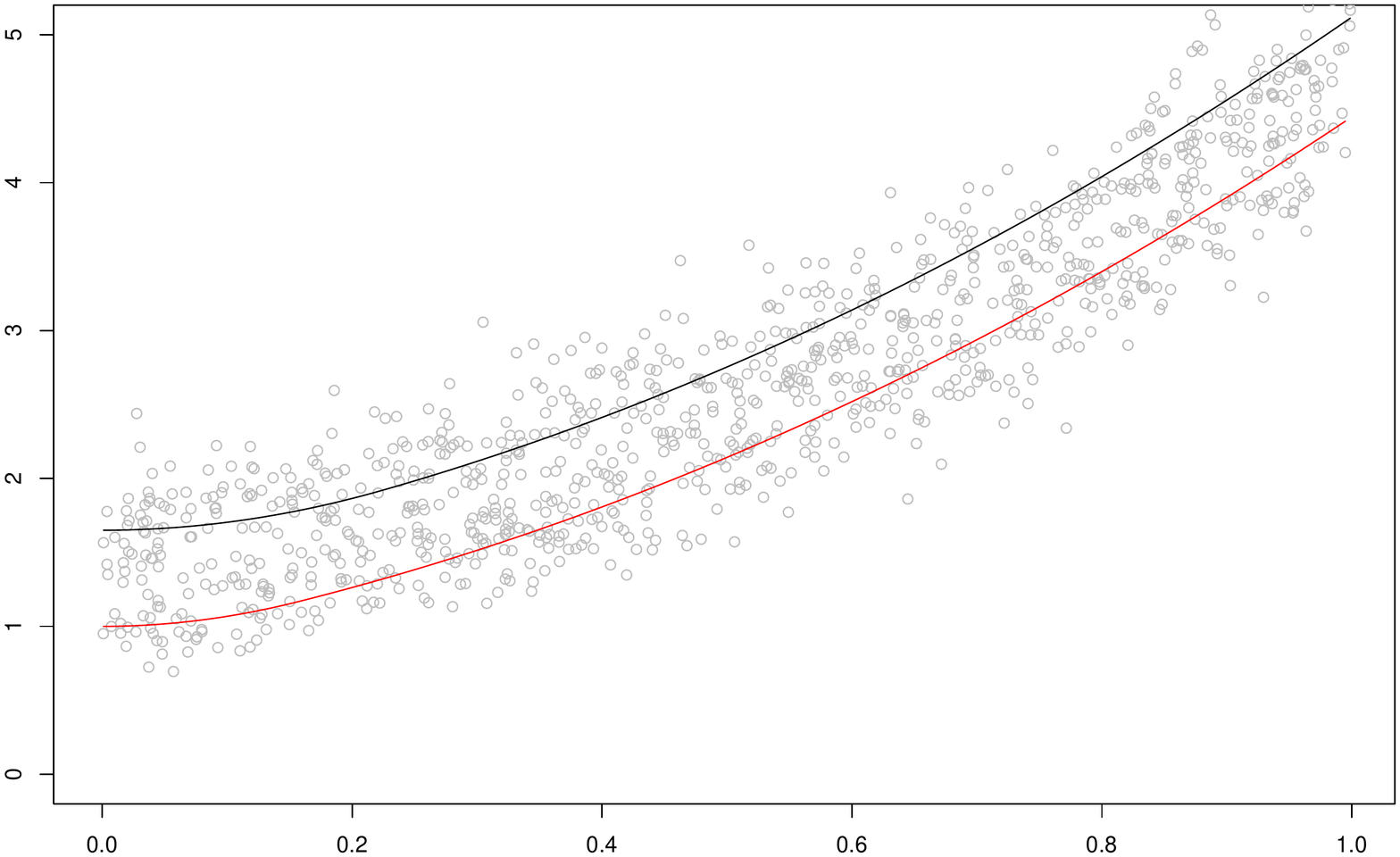}
\includegraphics[height=8cm,width=15cm,angle=0]{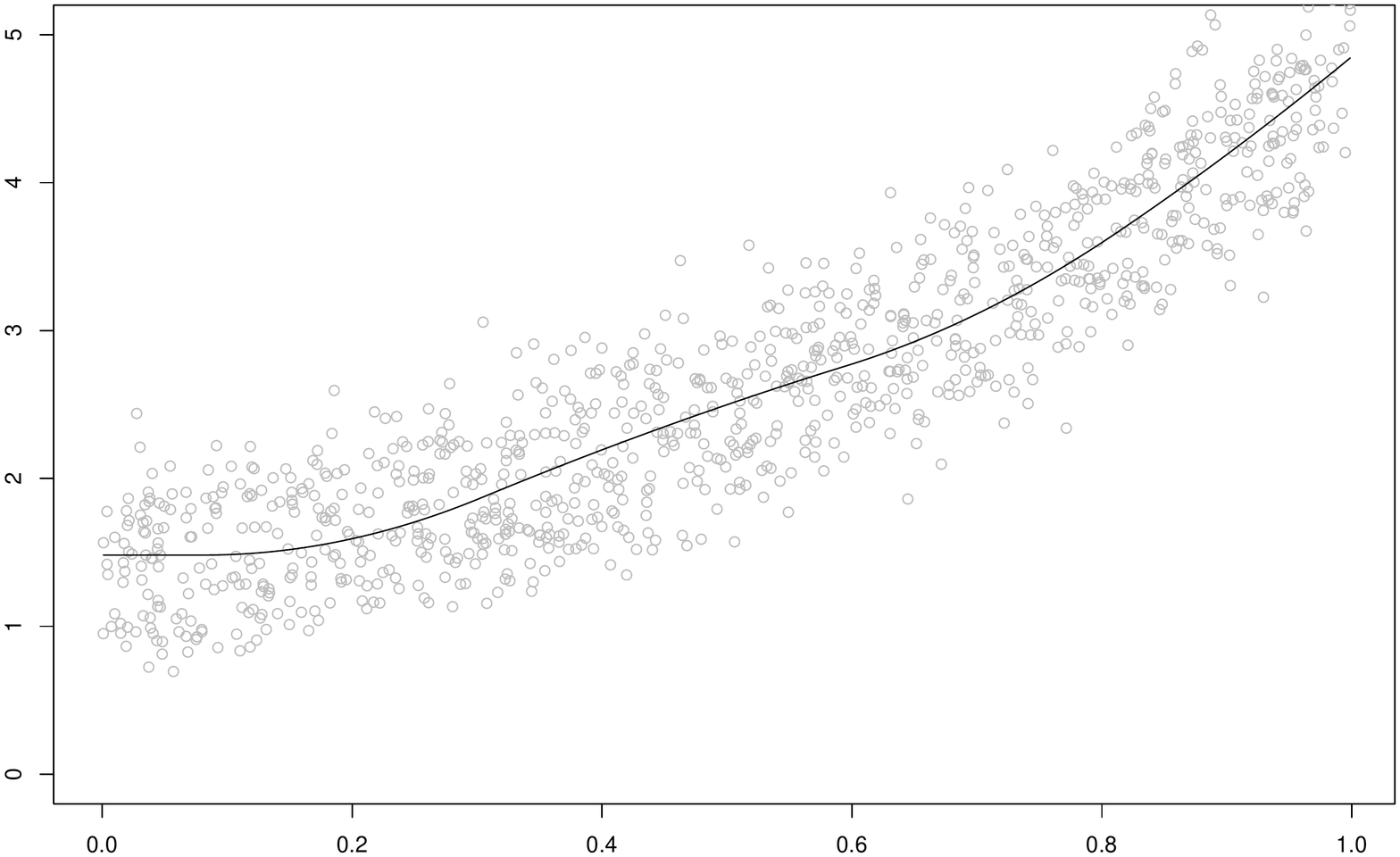}
\end{center}
\caption{The upper panel shows the results of minimizing the total
  variation of the second derivative for the data of the upper panel
  of Figure \ref{fig2} subject to monotonicity. The lower panel shows
  the corresponding result for the lower panel of Figure
  \ref{fig2}. \label{fig2a}}  
\end{figure}

\begin{figure}[hb]
\begin{center}
\includegraphics[height=8cm,width=15cm,angle=0]{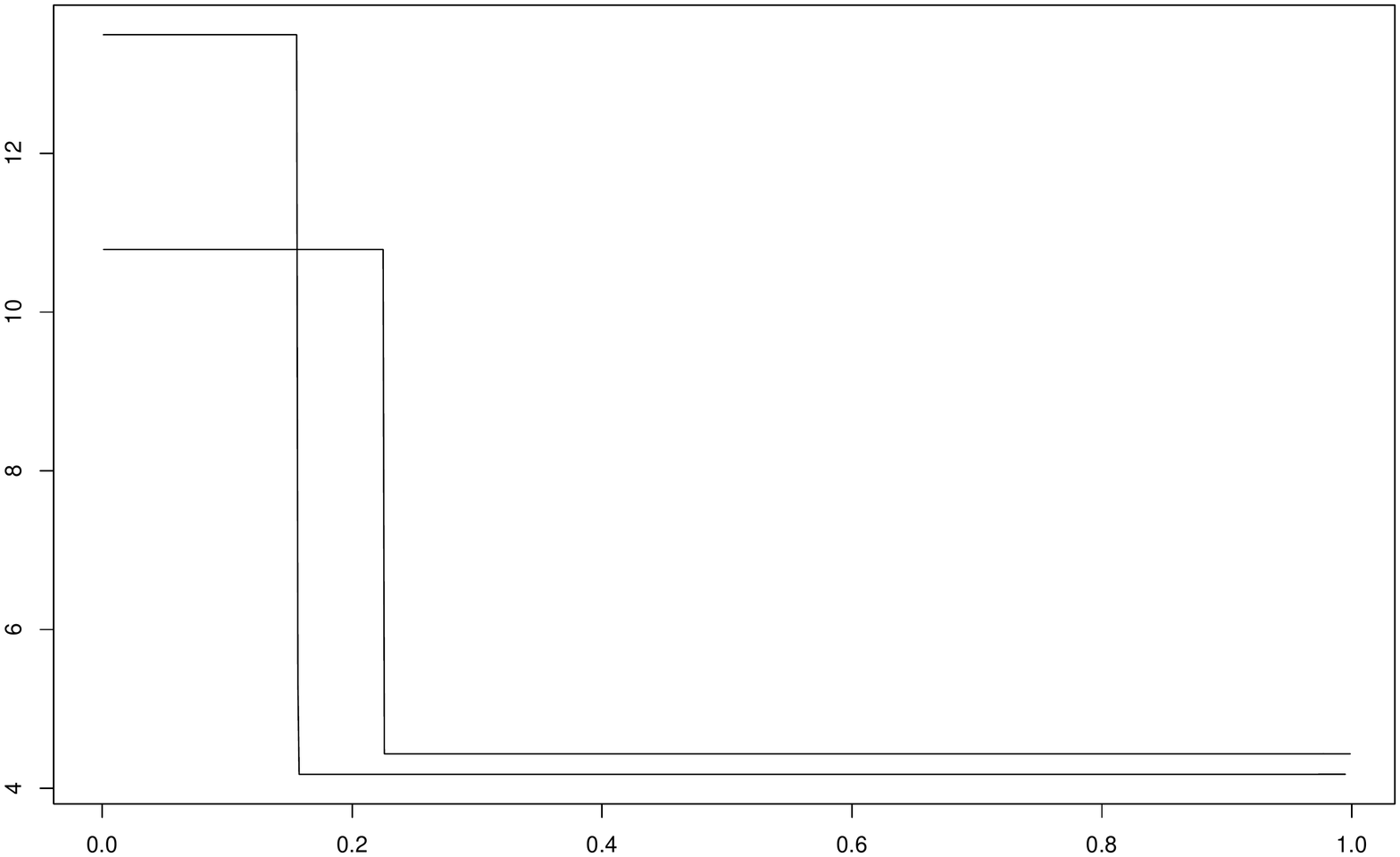}
\includegraphics[height=8cm,width=15cm,angle=0]{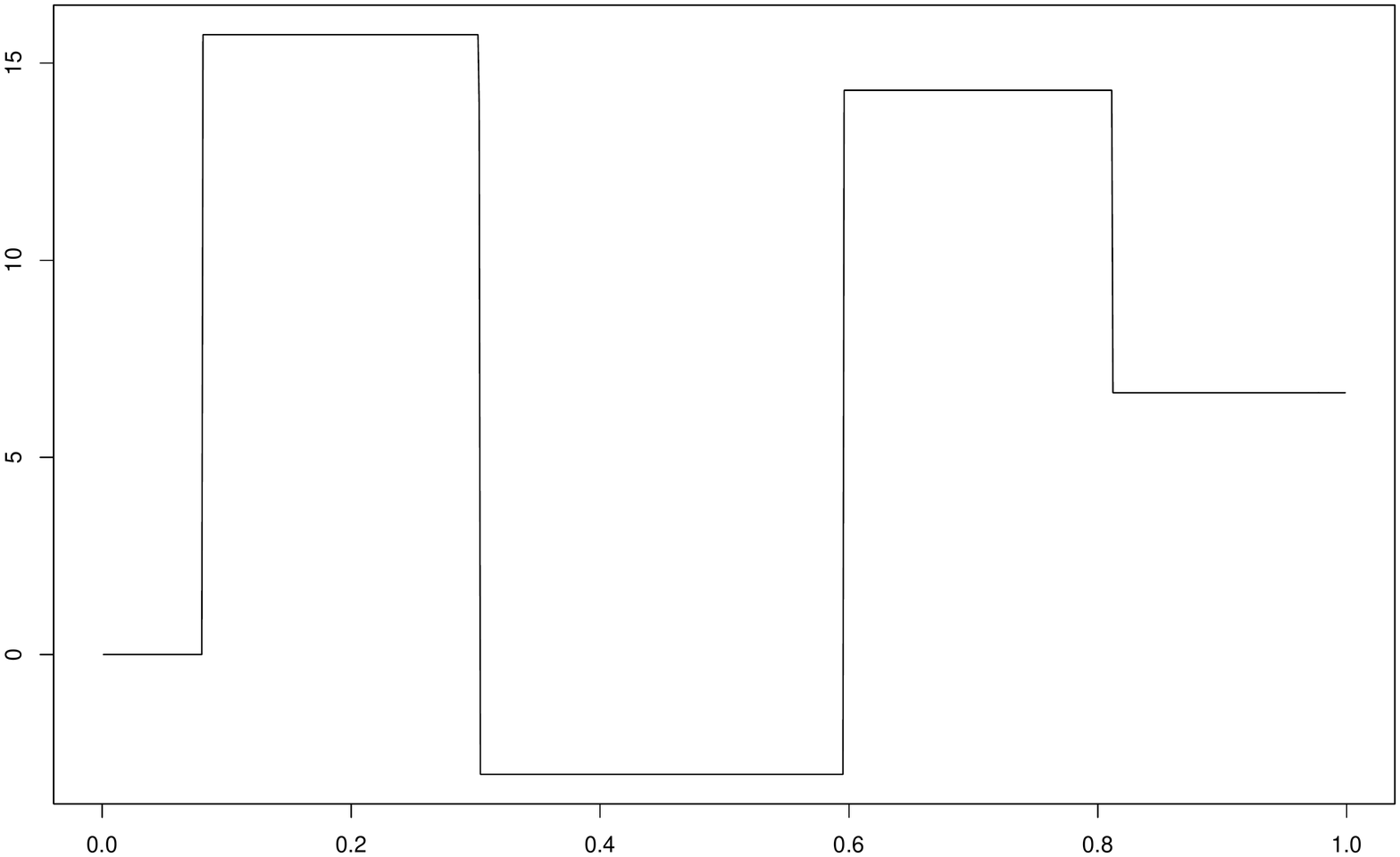}
\end{center}
\caption{The upper panel shows the second derivative iof the functions
  in the upper panel of Figure \ref{fig2a}: the lower panel shows
  the corresponding result for the lower panel of Figure
  \ref{fig2a}. \label{fig2b}}  
\end{figure}

Alternatively, rather than impose prior restrictions, one can determine a
simplest function in the joint approximating region. One possibility is
to minimize the number of local extreme points of a function $g$ subject
to $g \in {\mathcal A}_{{\mathbold n}_k}.$  Figure \ref{fig1} shows an
example of this approach where the joint approximating function has
514 internal local extreme points compared with monotone approximating
functions for both data sets separately. The additional local extreme
points can be regarded as the cost of the joint approximation. The
same idea can be used if simplicity is defined in terms of smoothness,
for example by minimizing the total variation $TV(g^{(2)})$  of the
second derivative subject to $g$ lying in the approximation
region. The upper panel of Figure \ref{fig2} shows the data and curves
of Figure \ref{fig1} but with the values of the second sample reduced
by an amount 2.3. There is now a joint monotone approximation which is
shown in the lower panel of Figure \ref{fig1} so there is no cost in
terms of the number of local extreme values. If we minimize the total
variation of the second derivative subject to the function being an
adequate monotone approximation then there is a cost. The upper panel
of Figure \ref{fig2a} shows the approximations of the two individual
samples  which minimize the total variation of the second derivatives
subject to the functions being monotone. The lower panel of Figure
\ref{fig2a} shows the joint approximation for the combined sample. The
second derivatives are shown in Figure \ref{fig2b}. The values of the
total variation of the second derivative are are he values  
are 9.317 and 6.305 for the individual samples and 59.496 for the
joint sample.

\subsection{Intersecting supports}
\begin{figure}[hb]
\begin{center}
\includegraphics[height=8cm,width=15cm,angle=0]{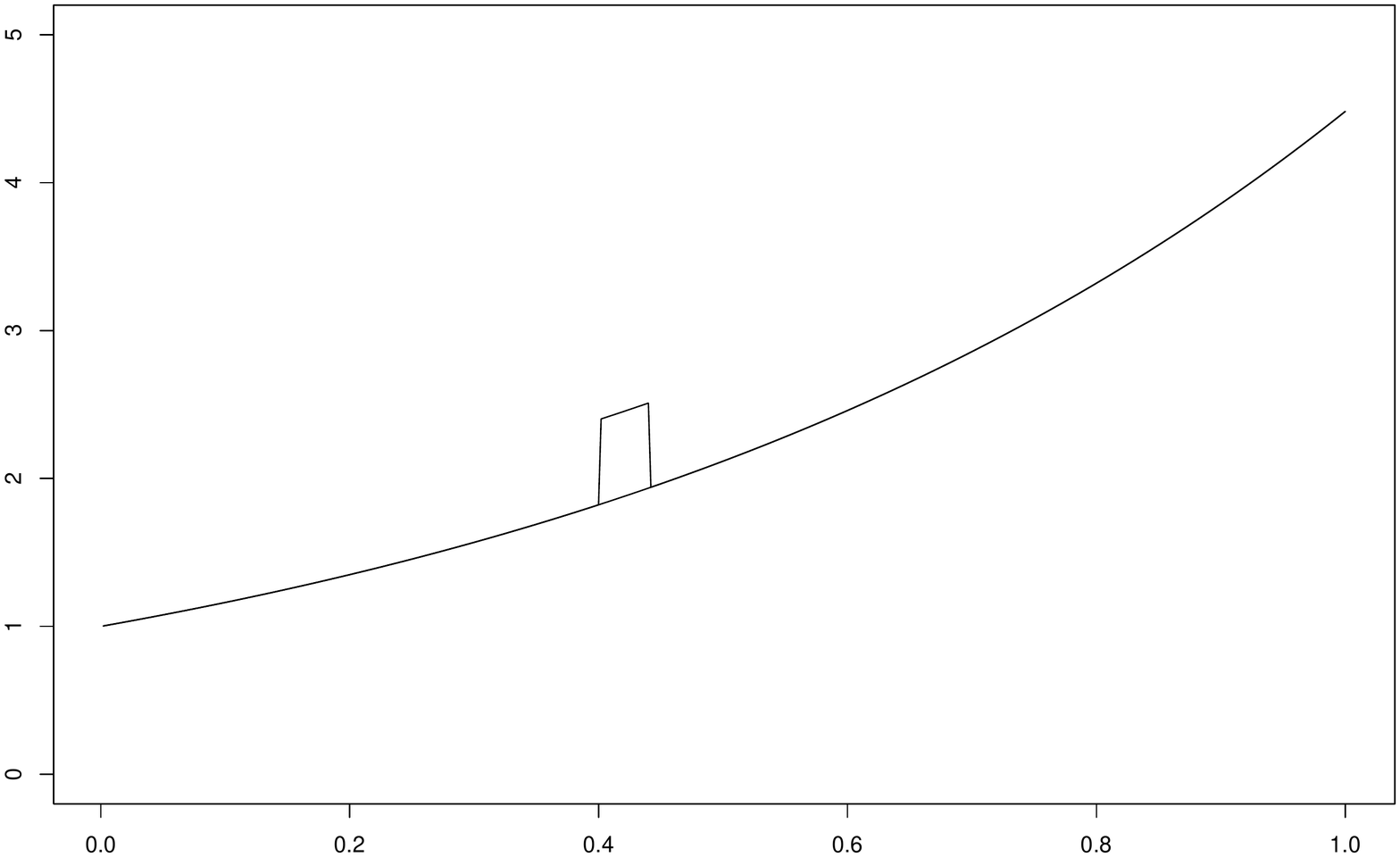}
\includegraphics[height=8cm,width=15cm,angle=0]{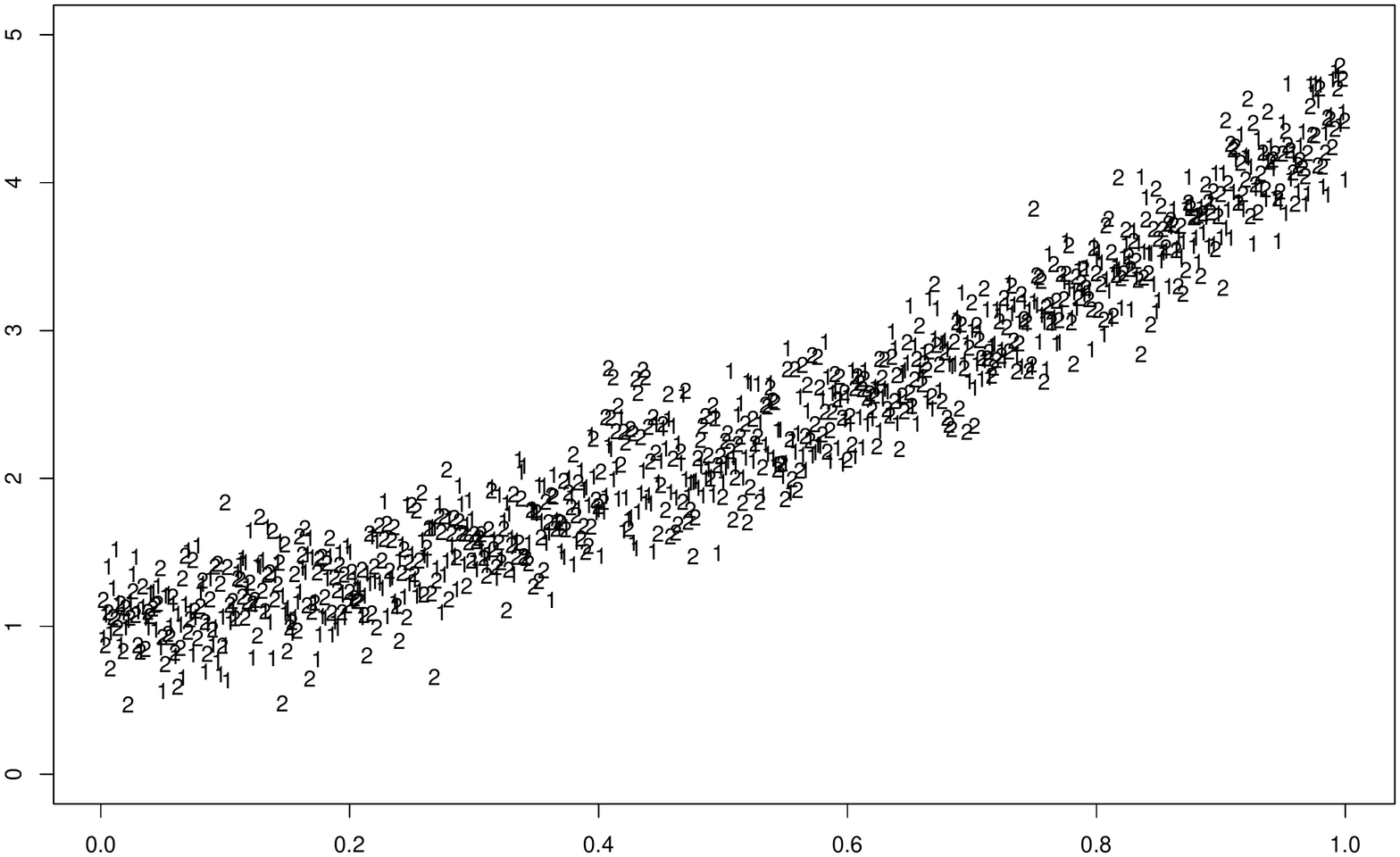}
\end{center}
\caption{The upper panel shows the function $f_1(t)=\exp(1.5t)$ and
  the function $f_2$ which is equal to $f_1$ apart from the interval
  $[0.402,0.440]$ where $f_2(t)=f_1(t)+0.575$.
 The lower panel shows the two data sets $Y_1(t_j)=f_1(t_j)+0.25Z_1(t_j)$
 and $Y_2(t_j)=f_2(t_j)+0.25Z_2(t_j)$ for $j=1,\ldots,500$ and with
 $t_j=j/500.$ \label{fig3}} 
\end{figure}

As mentioned in Section  \ref{introduction} the Neumeyer and Dette (2003)
procedure can detect differences of the order of $n^{-1/2}.$ We now
consider the size of detectable differences for our procedure in the
case of equal supports. For simplicity we consider only the case $k=2$
and assume that the  supports $S_{1n_1}$ and $S_{2n_2}$ are given by
$t_{1i}=t_{2i}=i/n.$ We take ${\mathcal I}_n$ to be
the set of all intervals but indicate below the adjustments required if
${\mathcal I}_n={\mathcal I}_n(\lambda)$ as in (\ref{multischeme}). We
state the results using $\sigma_1$ and $\sigma_2$ rather than the
estimates (\ref{sigma}) and write $\tau_n=\tau_n(\alpha^{1/k}).$ If a
joint approximating function ${\tilde f}_n$ exists then for any
interval $I$ of $[0,\,1]$ we have   
\begin{displaymath}
\frac{1}{\sqrt{\vert I\vert}}\left \vert \sum_{t_i \in I}
  (Y_j(t_i)-f_n(t_i))\right \vert \le \sigma_i\sqrt{\tau_n
  \log(n)}\,, \,\, j=1,\,2.
\end{displaymath}
and hence
\begin{displaymath}
\frac{1}{\sqrt{\vert I\vert}}\left \vert \sum_{t_i \in I}
  (Y_1(t_i)-Y_2(t_i))\right \vert \le
(\sigma_1+\sigma_2)\sqrt{\tau_n\log(n)}\,\,.
\end{displaymath}
For the noise we have with probability $\alpha$ 
\begin{displaymath}
\frac{1}{\sqrt{\vert I\vert}}\left \vert \sum_{t_i \in I}
  (Z_1(t_i)-Z_2(t_i))\right \vert \le (\sigma_1+\sigma_2)\,\sqrt{\tau_n
  \log(n)}
\end{displaymath}
and hence with probability $\alpha$ 
\begin{displaymath}
\frac{1}{\sqrt{\vert I\vert}}\left \vert \sum_{t_i \in I}
  (f_1(t_i)-f_2(t_i))\right \vert \le 2(\sigma_1+\sigma_2)\sqrt{\tau_n
  \log(n)}\,\,.
\end{displaymath}
Suppose now that $f_1$ and $f_2$ differ by an amount $\eta_n$ on an
interval $I_n \subset [0,\,1]$, that is $f_1(t)-f_2(t) > \eta_n, t\in
I_n$ and that the length of $I_n$ is $\delta_n.$ As $I_n$ contains about
$n\delta_n$ support points we see that 
\begin{displaymath} 
\frac{1}{\sqrt{n\delta_n}} n\delta_n \eta_n \le
2(\sigma_1+\sigma_2)\sqrt{\tau_n \log(n)}
\end{displaymath}
which implies that no joint approximation will exist if
\begin{equation} \label{sampsize1}
\sqrt{\delta_n}\,\eta_n > 2(\sigma_1+\sigma_2)\sqrt{\tau_n \log(n)/n}.
\end{equation}
It follows that with probability of at least $\alpha$, deviations
satisfying the latter inequality will be detected. If ${\mathcal
  I}_n={\mathcal I}_n(\lambda)$ as in (\ref{multischeme}) it follows
that there exists an interval $I_n'\subset I$ in ${\mathcal
  I}_n(\lambda)$ and of length $\vert I_n'\vert
\ge \vert I_n\vert/\lambda=\delta_n/\lambda$ for which $f_1(t)-f_2(t)
>\eta_n,\, t \in I_n'.$ This requires replacing (\ref{sampsize1}) by
\begin{equation} \label{sampsize1lam}
\sqrt{\delta_n}\,\eta_n >
2\sqrt{\lambda}(\sigma_1+\sigma_2)\sqrt{\tau_n \log(n)/n}. 
\end{equation}
We consider a situation similar to that of Figure \ref{fig1} as is
shown in Figure \ref{fig3}. The sample sizes are $n=500$ with
common supports $t_j=j/n$ and we take $\alpha$ to be 0.95 so that
$\alpha_k=0.95^{1/2}=0.9747.$ For this choice of $\alpha$  and with
${\mathcal I}_n={\mathcal I}_n(2)$ simulations give $\tau_n=2.973.$ We
set $f_1(t)=\exp(1.5t)$ and put $f_2(t)=f_1(t)$ except for $t \in
[0.402,0.44]$ where $f_2(t)=f_1(t)+\eta_n.$ For this interval
$\delta_n=20/500$ and so we expect to be able to detect deviations
$\eta_n$ of the order  
\begin{equation} \label{dev}
\eta_n=2\sqrt{2}(0.25+0.25)\sqrt{2.973\log 500\,}/\sqrt{20\,}=1.359
\end{equation}
with probability of at least 0.95. For the data shown in Figure
\ref{fig3} the difference is detected with $\eta_n=0.575$ but not with
$\eta_n=0.574.$ 

If we put $\delta_n=1$ in (\ref{sampsize1}) so that the two functions
deviate over the whole interval then
\begin{equation} \label{sampsize2}
\eta_n > 2\sqrt{2}(\sigma_1+\sigma_2)\sqrt{\tau_n \log(n)/n\,}.
\end{equation}
which implies that deviations of order $\sqrt{\log(n)/n\,}$ can be detected. 

The same analysis can be carried through using the approximation
region ${\mathcal A}_n^*.$ Corresponding to (\ref{sampsize1}) we
obtain
\begin{equation} \label{sampsize3}
\eta_n h(\sqrt{\delta_n})\, > 2(\sigma_1+\sigma_2)/\sqrt{n}
\end{equation}
where
\begin{equation}
h(\delta)=\frac{\delta}{2*\log(1/\delta)+
  \gamma_n(\alpha)\log\log(e^e/\delta)}\,, \,\, 0 <\delta\le 1,
\end{equation}
is monotonically increasing and $\gamma_n(\alpha)$ is bounded in $n$. In
particular, if $\delta_n=1,$ then deviations of order $1/\sqrt{n}$ can
be detected.

\subsection{Adapting the taut string algorithm}
The taut string algorithm of  \cite{DAVKOV01} has proved to be very
effective in determining the number of local extremes of a function
contaminated by noise (see  \cite{DAVGATNORWEI08}). We show how it
may be adapted to the case of $k$ samples. Let $n \le
\sum_{i=1}^kn_{i}$ denote the number of different $t_{ij}$ values
which we order as $0\le t_1<\ldots < t_n \le 1.$ For each sample we
calculate the values of ${\hat \sigma_{in_i}}$ given by (\ref{sigma}). We put 
\begin{equation} 
y_m=\frac{\sum_{t_{ij}=t_m}y_{ij}/{\hat \sigma_{in_i}}^2}{
  \sum_{t_{ij}=t_m}1/{\hat \sigma_{in_i}}^2}\,,\quad m=1,\ldots,n. 
\end{equation}
As a first step we check whether a joint approximating function
exists. We do this by putting ${\tilde f}_n(t_m)=y_m$ and then
determining whether ${\tilde f}_n \in {\mathcal A}_{{\mathbold n}_k}$. If
this is not the case we conclude that no joint approximation
exists. If a joint approximation exists we put  ${\hat
  \Sigma}_m=\sum_{t_{ij}=t_m}1/{\hat \sigma_{in_i}}^2$ and calculate
the partial sums $y^{\circ}_m=\sum_{1}^{m}{\hat
  \Sigma}_j y_j$ and  $A_m=\sum_{1}^{m}{\hat
  \Sigma}_j$ for $m=1,\ldots,\,n$ with $y^{\circ}_0=A_0=0$.
The initial lower and upper bounds $L_i$ and $U_i$ are set to be 
$L_i=Y_i-D,\, U_i=Y_i+D$ where $D$ is chosen so large that the
straight line joining $(0,0)$ and $(y^{\circ}_n,A_n)$ lies in the
tube. For a given tube, the taut string through the tube and constrained to
pass through $(0,0)$ and  $(y^{\circ}_n,A_n)$ is calculated. The value
of the estimate  ${\tilde f}_n(t_m)$ at the point $t_m$ is taken to be
the left hand derivative of the taut string except for the first point
where the right-hand derivative is taken. For each data set
individually it is now checked whether ${\tilde f}_n \in {\mathcal
  A}_{in_i}, i=1,\ldots,k.$ If this is the case the procedure
ends. Otherwise those intervals for which the inequalities defining
the ${\mathcal A}_{in_i}$ do not hold are noted and the tube is
squeezed at all points $t_{j-1}$ and $t_j$ for which $t_j$ lies in
such an interval. This is continued until a function ${\tilde f}_n \in
{\mathcal A}_{{\mathbold n}_k}$ is found.

\section{Comparison with other procedures} \label{comproc} 
\subsection{Analysis and simulations}
As the approach developed in this paper is not comparable with others
when the supports are disjoint, we restrict attention to the case of
equal supports.  For simplicity we take $k=2.$ For such data Delgado
(1992) proposed the test statistic   
\begin{equation} \label{cumsum1}
T_n=\sqrt{n}\max_{1\le j \le n}\vert
R(j)\vert/s_n^*=\max_{1\le j \le n}\left\vert
\sum_{i=1}^j(Y_1(t_i)-Y_2(t_i))\right\vert/(\sigma_n\sqrt{n}) 
\end{equation}
where $\sigma_n$ is some quantifier of the noise. Under the null
hypothesis $f_1=f_2=f$ the distribution of $T_n$ does not depend on
$f$. In this special case the test statistic of Neumeyer and Dette
also reduces to (\ref{cumsum1}). If the data were generated under
(\ref{basicmod}) then under $H_0$ the distribution of $T_n$
converges to that of $\max_{0 \le t \le 1}\vert B(t)\vert$
where $B$ is a standard Brownian motion. The 0.95-quantile is
approximately 2.24 which leads to rejection of $H_0$ if
\begin{equation} \label{cumsum2}
T_n \ge 2.24.
\end{equation}
Suppose now that the data are generated as in (\ref{basicmod}) with
$f_1(t)=f_2(t)$ apart from $t$ in an interval $I$ of length $\delta_n$
where $f_1(t)-f_2(t) \ge \eta_n$. It follows from (\ref{cumsum2})
that $H_0$  will be rejected with high probability if
\begin{equation} \label{lstdctdiff}
\delta_n\eta_n \ge 4.48\sigma/\sqrt{n}  
\end{equation}
where $\sigma^2=\sigma_1^2+\sigma_2^2.$ If $\delta_n=1$  deviations of
the order of $\sigma/\sqrt{n}$ can be picked up which contrasts with
the $O(\sigma\sqrt{\log(n)/n})$ of (\ref{sampsize2}).  If however
$\delta_n= 1/\sqrt{n}$ it follows from (\ref{lstdctdiff}) that the
test statistic $T_n$ will pick up deviations of the order of
$\sigma$. It follows from (\ref{sampsize1}) and (\ref{sampsize3}) that
the methods based on  ${\mathcal A}_n$ and ${\mathcal A}_n^*$ will
both pick up deviations of the order of $\sigma\sqrt{\log(n)/\sqrt{n}\,}.$

Another test which is applicable in this situation is due to Fan and 
Lin (1998). If we denote the Fourier transform of the data sets by
${\tilde Y}_1(i)$ and ${\tilde Y}_2(i), i=1,\ldots, n,$ and order
them  as described in Fan and Lin (1998), their test statistic reduces
to 
\begin{equation} \label{fanlint}
T^*_n =\max_{1 \le m \le n} \left\vert\frac{1} {\sqrt{m}}\sum_{i=1}^m
  (({\tilde Y}_2(i)-{\tilde Y}_1(i))^2/{\tilde \sigma}_n^2-1)\right\vert  
\end{equation}
where ${\tilde \sigma}_n$ is some estimate of the standard deviation
of the ${\tilde Y}_2(i)-{\tilde Y}_1(i).$  For data generated under the model
(\ref{basicmod}) the critical value of $T^*_n$ can be obtained by
simulations. It is not as simple to determine the size of the
deviations which can be detected by the test (\ref{fanlint}) as the
test statistic is a function of the Fourier transforms and the
differences in the functions must be translated into differences in
the Fourier transforms. The first member of the sum in (\ref{fanlint})
is the difference of the means and this is given the largest
weight. We do not pursue this further but give the results of a small
simulation study. 

\begin{figure}[hb]
\begin{center}
\includegraphics[height=5cm,width=9cm,angle=0]{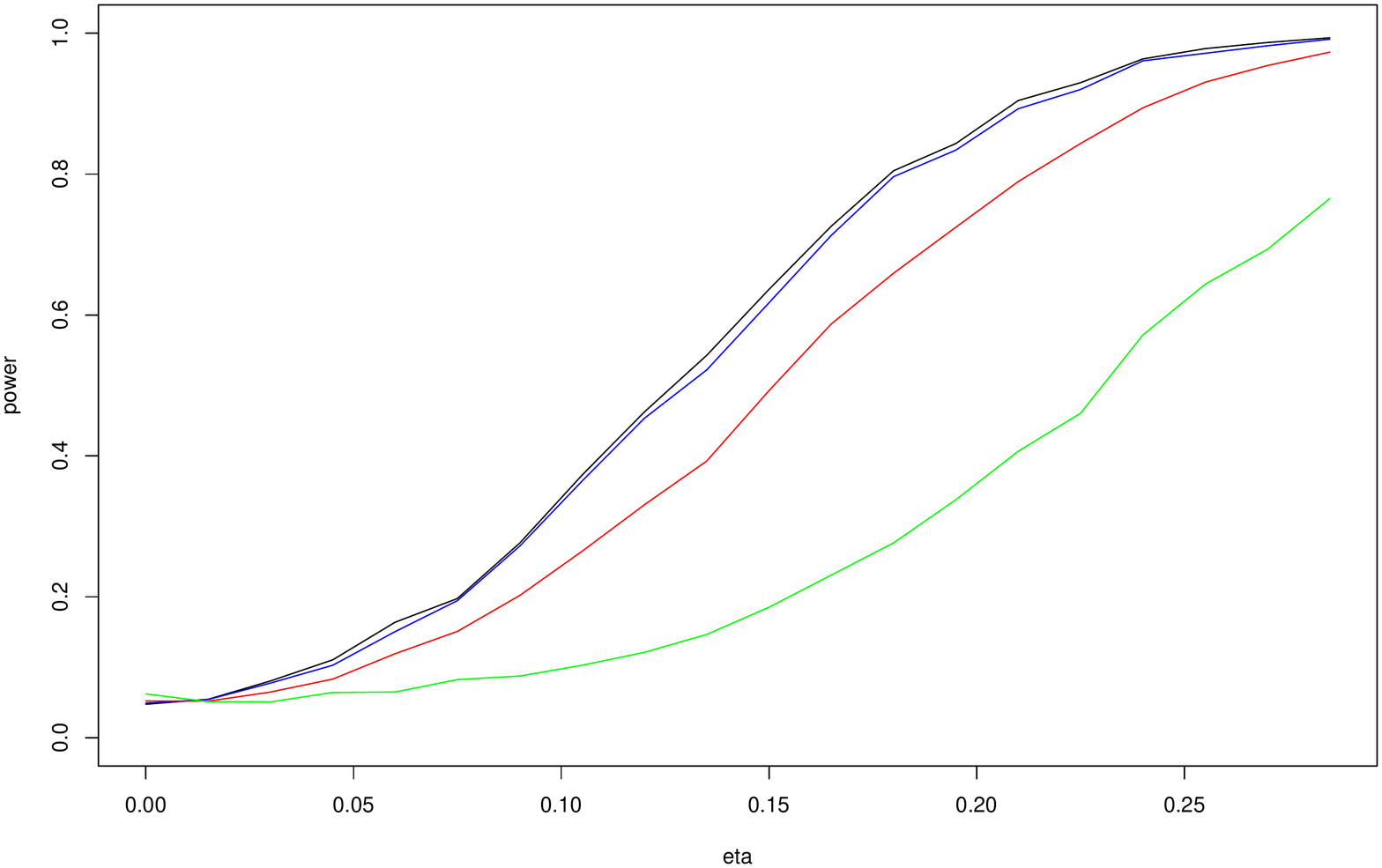}
\includegraphics[height=9cm,width=5cm,angle=90]{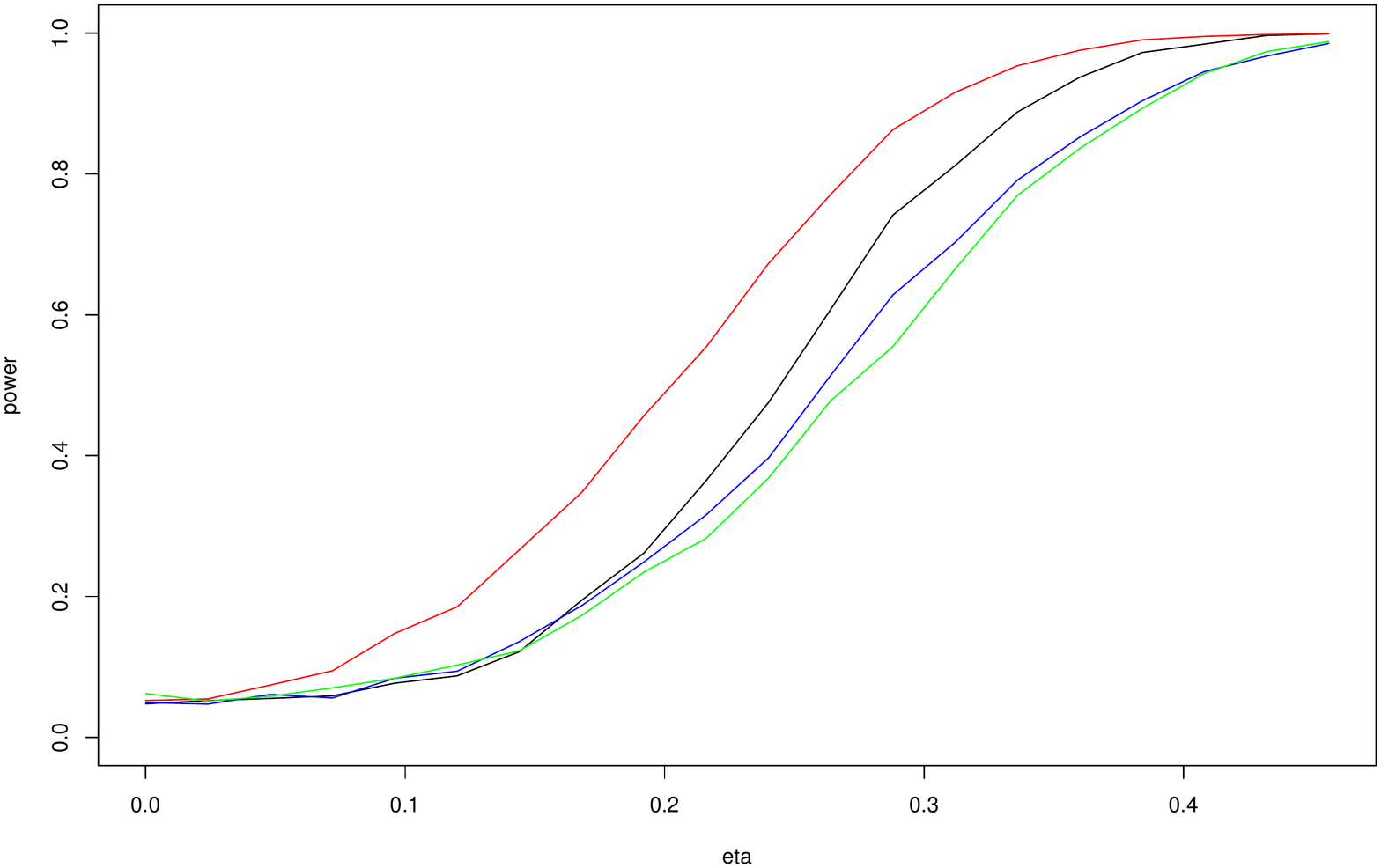}\\
\includegraphics[height=9cm,width=5cm,angle=90]{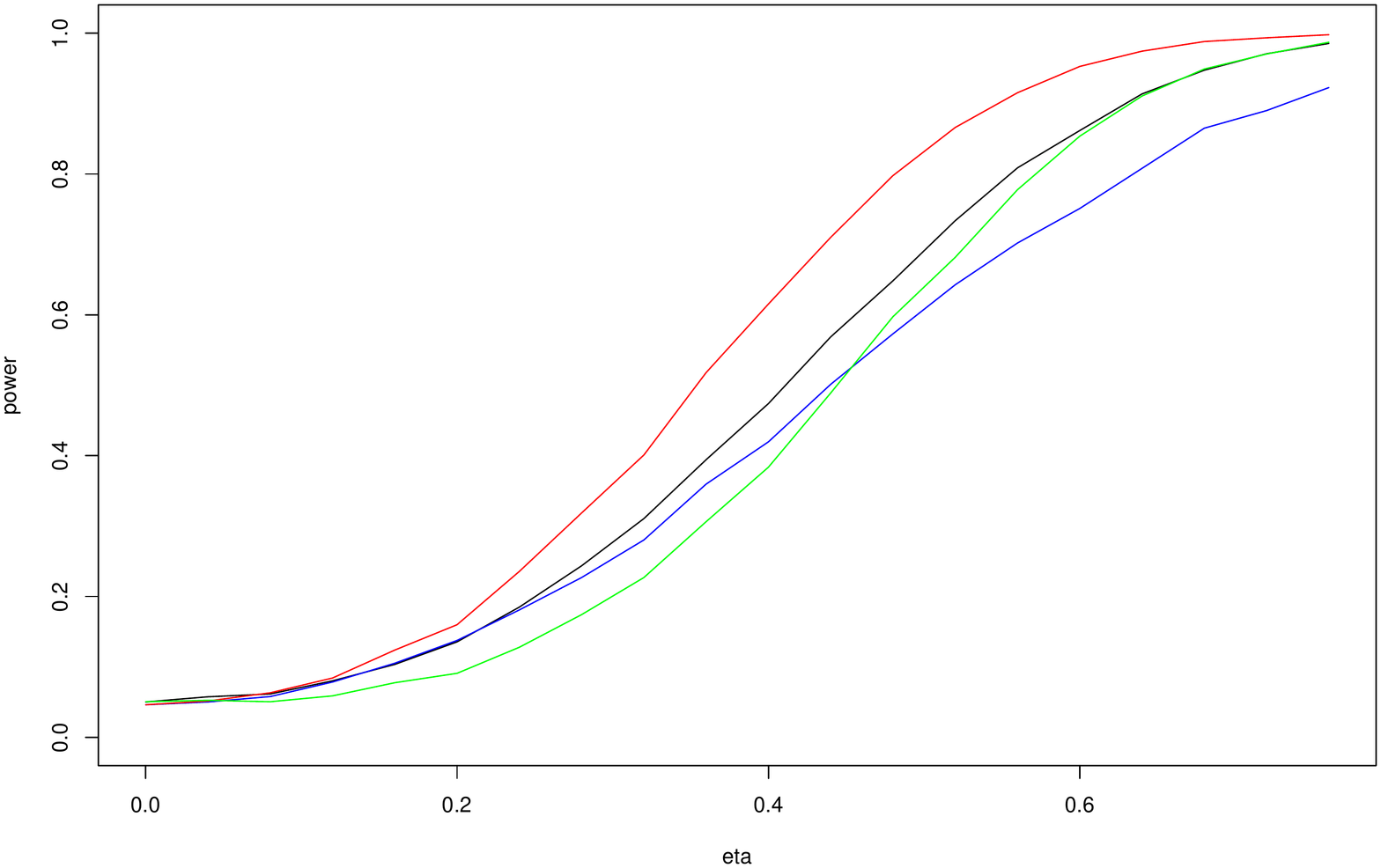}
\includegraphics[height=9cm,width=5cm,angle=90]{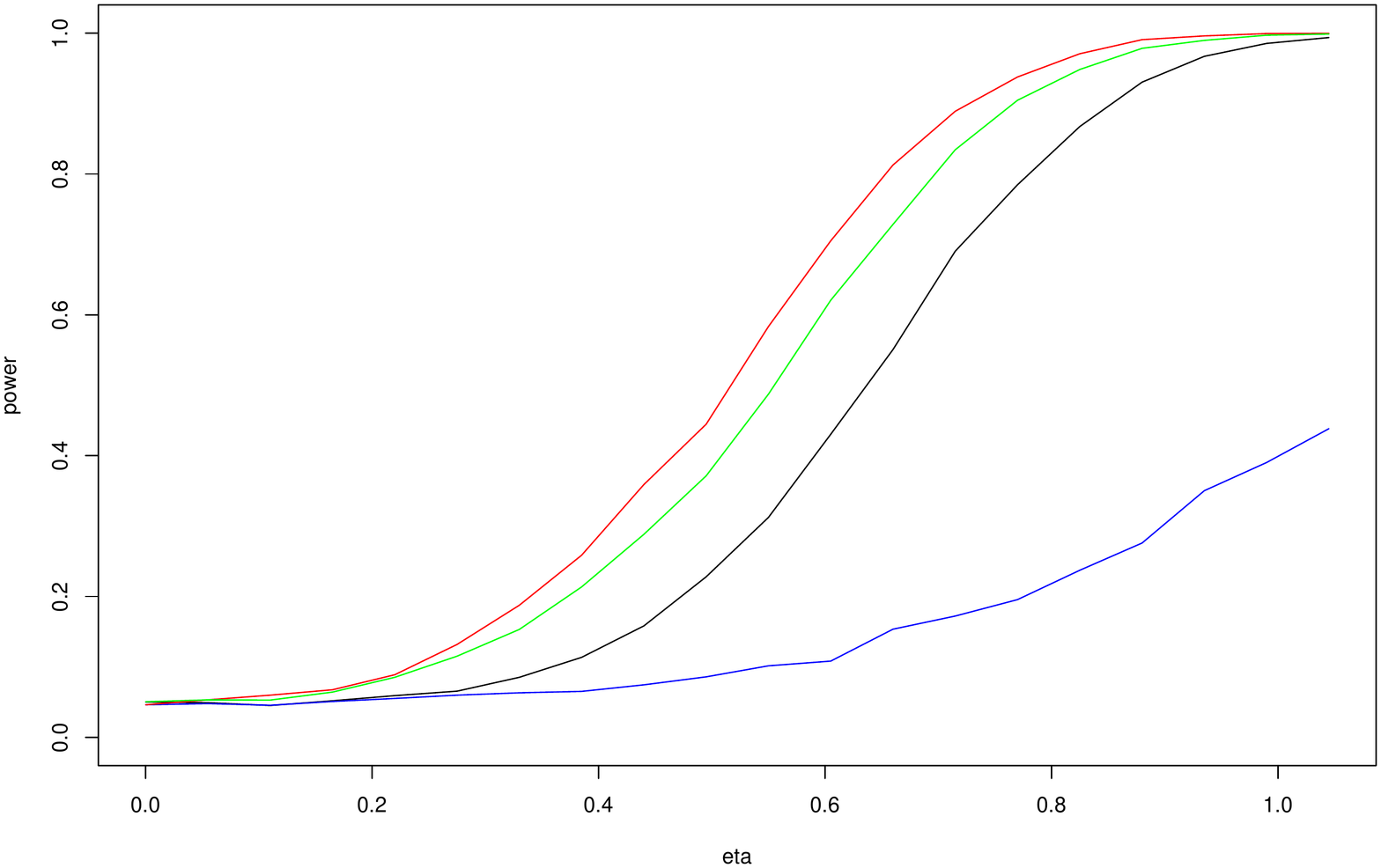}\\
\end{center}
\caption{The top row shows the power functions of the four tests with
  $g$ given by (1) and (2). The bottom row shows the power functions
  of the four tests with   $g$ given by (1) and (2). The
  Delgado--Neumeyer--Dette  is shown in blue, the Fan--Lin test in
  black, the test based on ${\mathcal A}_n$ in green and that based on
  ${\mathcal A}_n^*$ in red. \label{power}} 
\end{figure}

We put $n=500$ and consider two samples of the form
\begin{eqnarray}
Y_1(i/n)&=& Z_1(i/n),\quad i=1,\ldots, n=500 \label{ysamp1}\\
Y_2(i/n)&=& g(i/n)+Z_2(i/n), \quad i=1,\ldots, n=500 \label{ysamp2}
\end{eqnarray}
 were generated  where the $Z_j(i/n)$ are i.i.d $N(0,1)$ random
 variables with $g$ given by one of
\begin{equation} \label{varg}
\begin{tabular}{ll}
(1)\quad $g_1(t)=\eta,\, 0 \le t \le 1,$& (2)\quad $g_2(t)=\left \{
  \begin{array}{ll} \eta,&0 \le t \le 1/2,\\
-\eta. &1/2 < t \le 1,\\ 
\end{array} \right.$\\ 
&\\
(3)\quad $g_3(t)=\left \{
  \begin{array}{lll}
0,& 0 \le t \le U,\\ 
\eta,&U < t \le U+1/4,\\
0,& U+1/4 < t \le 1,\\
\end{array}\right.$
&(4)\quad $g_4(t)=\left \{
  \begin{array}{lll}
0,& 0 \le t \le U,\\ 
\eta,&U < t \le U+1/8,\\
-\eta,&U +1/8< t \le U+1/4,\\
0,& U+1/4 < t \le 1,\\
\end{array}\right.$\\
&\\
\end{tabular}
\end{equation}
where $U$ is uniformly distributed on $[0,3/4]$ and independent of the
$Z_i, i=1,2.$ The four procedures, Delgado--Neumeyer--Dette, Fan--Lin
and those based on ${\mathcal A}_n$ and  ${\mathcal A}_n^*$ were
all calibrated to give tests of size $0.05$ for testing $g\equiv 0.$
The critical values for  Delgado--Neumeyer--Dette and Fan--Lin tests
are 2.22 and 6.97 respectively. The value of $\tau_n$ for the test
based on ${\mathcal A}_n$ is 1.46 and the corresponding value of
$\gamma_n$ for that based on ${\mathcal A}_n^*$ is 0.66.
Figure \ref{power} shows the power of the tests for different values
of $\eta.$ The upper panels are the results for $g$ given by
(\ref{varg}) (1) and (\ref{varg}) (2) and the
lower panels for $g$ given by (\ref{varg}) (3) and (\ref{varg})
(4). The colour scheme is as follows: Delgado--Neumeyer--Dette blue,
the Fan--Lin black, ${\mathcal A}_n$ green and  ${\mathcal
  A}_n^*$ red. The results confirm the analysis given above. The
Delgado--Neumeyer--Dette and Fan--Lin tests are better with $g$ given
by (1) but if the mean difference is zero (2), or the interval is
small (3) or both (4) then they are outperformed by the procedure
based on ${\mathcal A}_n^*$ and, in case 4, also by that based on
${\mathcal A}_n.$ 

\subsection{An application}
\begin{figure}[ht]
\begin{center}
\includegraphics[height=15cm,width=15cm,angle=0]{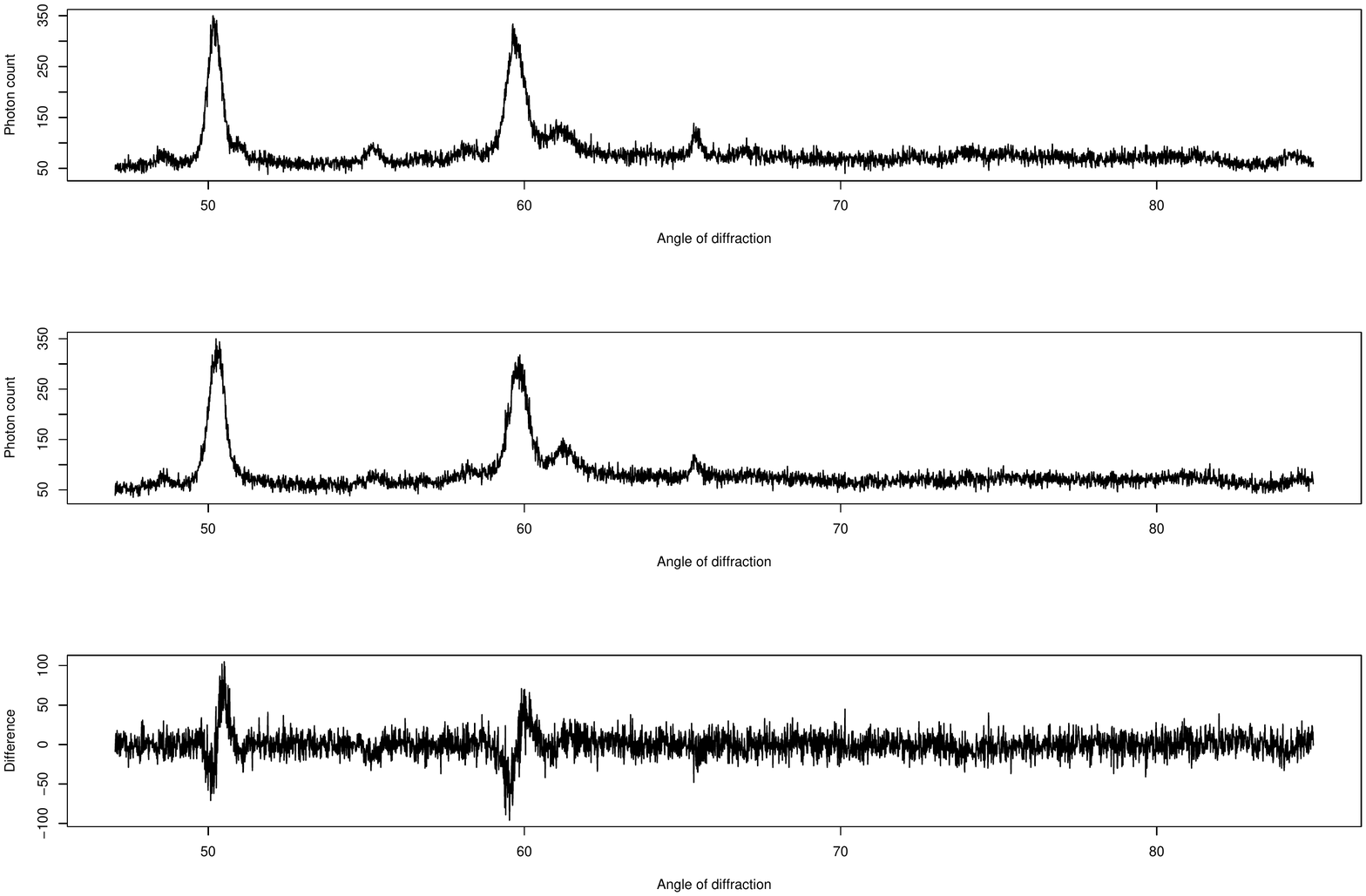}
\end{center}
\caption{The top and center panels show two data sets each of 4806
  observations with the same design points. The lower panel shows the
  differences of the two samples.\label{mergel1}}
\end{figure}

We give an example with some real data from the area of thin-film
physics. They give the number of photons of refracted X-rays
as a function of the angle of refraction and were kindly supplied by
Professor Dieter Mergel of the University of Duisburg-Essen.  Two such
data sets are shown in the top panel of Figure  \ref{mergel1}; the
differences $y_1(t_i)-y_2(t_i)$ are shown in the bottom panel. Each data
set is composed of 4806 measurements and the design points are the
same. The samples differ in the manner in which the thin film was
prepared. One of the questions to be answered is whether the results
of the two methods are substantially different.\\

The noise levels for the data sets are the same, namely 8.317, which
is explainable by the fact that the data are integer valued. The
differences between the two data sets are concentrated on intervals
each containing about 40 observations. The estimate (\ref{lstdctdiff})
suggests that the differences will have to be of the order of 92 to be
detected with a degree of certainty by the Delgado--Neumeyer--Dette
test. The actual differences are of about this order and in fact the
test fails to reject the null hypothesis at the 0.1 level. The
realized value of the test statistic is 1.734 as 
against the critical value of 1.90 given in (\ref{cumsum2}). The
Fan-Lin test (\ref{fanlint}) rejects the  null hypothesis at the 0.01
level. The realized value of the test statistic is 111.66 as against
the critical value of 12.44 for a test of size $\alpha=0.01.$ Finally
the tests based on ${\mathcal A}_n$ and ${\mathcal A}_n^*$ both reject
the null hypothesis at the 0.01 level. The realized value of
$\tau_n$ is 43.15 as against the critical value of 1.50. The realized
value of $\gamma_n$ is 53.27 as against the critical value of 0.733.  



\bibliographystyle{apalike}
\bibliography{npcomp}

\newpage

\newpage

\newpage
\newpage

\newpage

\newpage

\newpage

\end{document}